\documentclass[11pt]{amsart}
\usepackage{amsmath,paralist}
\usepackage{amsthm}
\usepackage{graphicx}
\usepackage{tikz,color}
\usepackage[margin=1in]{geometry}
\usetikzlibrary{arrows,shapes,snakes,backgrounds} 

\usepackage[all]{xy}

\newtheorem{theorem}{Theorem}
\newtheorem*{corollary*}{Corollary}
\newtheorem*{theorem*}{Theorem}

\newtheorem{proposition}[theorem]{Proposition}
\newtheorem{lemma}[theorem]{Lemma}
\newtheorem{corollary}[theorem]{Corollary}

\theoremstyle{remark}

\theoremstyle{definition}

\newcommand{\dd} {n}
\newcommand{\eq}{\begin{equation}}
\newcommand{\lb}[1]{\label{#1}}

\newcommand{\skk}{{{{\mbox{${\mbox{\scriptsize\boldmath$k$}}$}}}}}

\newcommand{\en}{\end{equation}}
\newcommand{\kk}{{\mbox{\boldmath$k$}}}
\newcommand{\nnints}{\mathbb{N}}

   \def\vol{{\rm vol}}
 
\def\v{{\rm v}}
\def\i{{\rm indeg }}
\def\ini{{\rm in }}
\def\fin{{\rm fin}}

 \def\f_H{{\bf w}}
 \def\f{{\bf f}}
 
\def\R{\mathbb{R}}

\def\Z{\mathbb{Z}}

 \def\F{\mathcal{F}}
 \def\P{\mathcal{P}}
\def\T{\mathcal{T}}
\def\is{{\rm InSeq}}
\begin{document}

\tikzstyle{w}=[label=right:$\textcolor{red}{\cdots}$] 
\tikzstyle{b}=[label=right:$\cdot\,\textcolor{red}{\cdot}\,\cdot$] 
\tikzstyle{bb}=[circle,draw=black!90,fill=black!100,thick,inner sep=1pt,minimum width=3pt] 
\tikzstyle{bb2}=[circle,draw=black!90,fill=black!100,thick,inner sep=1pt,minimum width=2pt] 
\tikzstyle{b2}=[label=right:$\cdots$] 
\tikzstyle{w2}=[]
\tikzstyle{vw}=[label=above:$\textcolor{red}{\vdots}$] 
\tikzstyle{vb}=[label=above:$\vdots$] 

\tikzstyle{level 1}=[level distance=3.5cm, sibling distance=3.5cm]
\tikzstyle{level 2}=[level distance=3.5cm, sibling distance=2cm]

\tikzstyle{bag} = [text width=4em, text centered]
\tikzstyle{end} = [circle, minimum width=3pt,fill, inner sep=0pt]

\title[Product formulas for  volumes of flow polytopes]{Product formulas for  volumes of flow polytopes}
\author{Karola M\'esz\'aros}
\address{
Department of Mathematics, University of Michigan, Ann Arbor, MI 48109
}
\thanks{The author is supported by a National Science Foundation Postdoctoral Research Fellowship}
\date{November 22, 2011}

\begin{abstract} 
Intrigued by the product formula $\prod_{i=1}^{n-2} C_i$ for the volume of the Chan-Robbins-Yuen polytope $CRY_n$, where $C_i$ is the $i^{th}$ Catalan number, we construct a family of polytopes $\mathcal{P}_{m,n}$, indexed by $m \in \mathbb{Z}_{\geq 0}$ and $n \in \mathbb{Z}_{\geq 2}$, whose volumes are given by the product   $$\prod_{i=m+1}^{m+n-2}\frac{1}{2i+1}{{m+n+i} \choose {2i}}.$$  The Chan-Robbins-Yuen polytope $CRY_n$ coincides with  $\P_{0,n-1}$. Our construction of the polytopes $\mathcal{P}_{m,n}$ is an application of a systematic method we develop for expressing volumes of a class of flow polytopes as the  number of certain triangular arrays.  This method can also be used as a heuristic technique for constructing polytopes with combinatorial volumes. As an illustration of this we construct polytopes whose volumes equal the number of $r$-ary trees on $n$ internal nodes, $ \frac{1}{(r-1)n+1} {{rn} \choose n}$.
Using triangular arrays we also express the volumes of  flow polytopes as  constant terms of formal Laurent series.  \end{abstract}

\maketitle
  
  \section{Introduction}
\label{sec:in}

In this paper we device an encoding of  triangulations for a large class of flow polytopes. Using this encoding,  we  prove new volume formulas for a family of polytopes, of which, the famous  Chan-Robbins-Yuen polytope is a special case. We also use the encoding as a heuristic for constructing polytopes with combinatorial volumes, as well as to shed light on the geometric nature of several intriguing conjectures of Chan, Robbins, and Yuen \cite{cry}, which previously had no polytopal interpretations. We tie the whole story in with  constant terms of formal Laurent series.

The Chan-Robbins-Yuen polytope $CRY_{n+2}$ was discovered in 1998  by Chan, Robbins, and Yuen \cite{cry},  while they were studying the Birkhoff polytope. (The $CRY_{n+2}$  polytope is a face of the Brikhoff polytope.) Chan, Robbins, and Yuen \cite{cry}  conjectured that the volume of $CRY_{n+2}$ was $\prod_{i=1}^n C_i$,  where $C_i=\frac{1}{2i+1}{2i \choose i}$ is the $i^{th}$ Catalan number. This conjecture, which thanks to Zeilberger  \cite{z} quickly  became a theorem,  has captivated combinatorialists and noncombinatorialists alike.  It continues to  intrigue, since no combinatorial proof of it is known.  

Without reference to polytopes the above well-known result  can be stated as 
\begin{equation} \label{a}K_{A_n^+}(1, 2, \ldots, n, -{n+1 \choose 2})=\prod_{k=1}^{n}C_k,\end{equation}
where $K_{A_n^+}(\v)$ is the Kostant partition function of type $A_n$, which is equal to the 
 number of ways to write the vector $\v$ as a nonnegative linear combination of the positive type $A_n$ roots  without regard to order. 
 
 In their paper \cite{cry} Chan, Robbins and Yuen defined the $CRY_{n+2}$ polytope and constructed a triangulation of it, the simplices of which they bijected  with certain triangular arrays. Their triangulation is specific to the $CRY_{n+2}$ polytope and provokes the question of how the triangular arrays came up. In their famous Conjecture 1, Chan, Robbins, and Yuen stated that the number of the triangular arrays encoding the triangulation of  $CRY_{n+2}$ is  $\prod_{i=1}^n C_i$. It was this form of the conjecture that Zeilberger \cite{z} proved  in a very elegant, though noncombinatorial manner, using the Morris constant term identity \cite{wm}.  Chan, Robbins and Yuen \cite[Conjecture 2, 3]{cry} gave two other conjectural formulas for the number of certain other triangular arrays, which they did not give  polytopal interpretations for and which they mentioned came up in their computational exprimentations. Zeilberger \cite{z} showed how to prove those conjectures, too, using  the Morris constant term identity \cite{wm}.  Among other things, in this paper we show that  the triangular arrays in Conjectures 2 and 3 of \cite{cry} can be interpreted in the world of polytopes. Furthermore, we show that  triangular arrays in general, of which the ones appearing in Conjectures 1,2 and 3 of \cite{cry} are just special cases in disguise, are an essential part of the study of triangulations of flow polytopes.

We develop a systematic method for expressing volumes of a class of flow polytopes as the  number of certain triangular arrays. As an application of this method we construct a family of polytopes $\mathcal{P}_{m,n}$, indexed by $m \in \mathbb{Z}_{\geq 0}$ and $n \in \mathbb{Z}_{>0}$, whose normalized  volumes are given by the product   $$\prod_{i=m+1}^{m+n-2}\frac{1}{2i+1}{{m+n+i} \choose {2i}}.$$    The Chan-Robbins-Yuen polytope $CRY_n$ coincides with  $\P_{0,n-1}$.    Since the volume of any flow polytope can be expressed as a Kostant partition function evaluated at a vector,  our results also imply interesting formulas for special evaluations of the Kostant partition function.  Among other identities, we deduce a generalization of \eqref{a} using the  polytopes $\mathcal{P}_{m,n}$ mentioned above:

\begin{equation} \label{d} K_{A_{n}}(m+1, m+2, \ldots, m+n, -nm-{n \choose 2})=\prod_{i=m+1}^{m+n-2}\frac{1}{2i+1}{{m+n+i} \choose {2i}}. \end{equation} 

Equation \eqref{d}  has been previously stated without proof  by  A.N. Kirillov in a  different, but equivalent form, in his beyond-amazing paper entitled ``Ubiquity of Kostka polynomials"  \cite{kir}.

Our  method of encoding triangulations as triangular arrays also  yields a heuristic technique for constructing polytopes with combinatorial volumes. As an illustration of this we construct polytopes whose volumes equal the number of $r$-ary trees on $n$ internal nodes, $ \frac{1}{(r-1)n+1} {{rn} \choose n}$. We also express the volumes of  flow polytopes as  constant terms of formal Laurent series.

The outline of this paper is as follows. In Section \ref{sec:back} we provide the basic definitions and  background about flow polytopes and Kostant partition functions. In Section \ref{sec:arr} we present an algorithm for triangulating flow polytopes and prove that the volume of a flow polytope is equal to the number of certain triangular arrays depending on the polytope. We also express the latter as the 
constant term of a formal Laurent series. Section \ref{sec:m} is devoted to studying a family of polytopes $\mathcal{P}_{m,n}$, $m \in \mathbb{Z}_{\geq 0}$, $n \in \mathbb{Z}_{\geq 2}$, 
of which the Chan-Robbins-Yuen polytope is a special case, and proving that the normalized volume of $\mathcal{P}_{m,n}$ is $$\prod_{i=m+1}^{m+n-2}\frac{1}{2i+1}{{m+n+i} \choose {2i}}.$$ At the same time we give \cite[Conjecture 3]{cry} a polytopal interpretation. In Section \ref{sec:conj2} we provide a polytopal interpretation of \cite[Conjecture 2]{cry}, thereby constructing collections of polytopes $\{P^{n,k}_i\}_{i \in C_{n,k}}$, $k \leq n$,  whose volumes sum to $$N(n,k)\times \prod_{i=1}^{n-1}C_i,$$ where $$N(n, k)=\frac{1}{n}{n \choose k}{n \choose {k-1}}$$ is the Narayana number.
    In Section  \ref{sec:constr} we construct a family of polytopes with volumes equal the number of $r$-ary trees on $n$ internal nodes, $ \frac{1}{(r-1)n+1} {{rn} \choose n}$, as an illustration of how to use our triangular array encoding for constructing polytopes with combinatorial volumes.   Section \ref{sec:cat} is meant as some pleasant little fun:  we study natural generalizations of  the polytope Catalonotope, and calculate their volumes. Postnikov and Stanley were the first to study (an equivalent version of)  the Catalanotope and  they calculated its volume,  the Catalan number; thus, its name. 
   
 \section{Background information}
 \label{sec:back}

\subsection{Flow polytopes.}Flow polytopes are associated to loopless graphs in the following way. Let $G$ be a graph on the vertex set $[n+1]$, and let $\ini(e)$ denote the smallest (initial) vertex of edge $e$ and $\fin(e)$ the biggest (final) vertex of edge $e$.  Think of fluid flowing on the edges of $G$ from the smaller to the bigger vertices, so that the total fluid volume entering vertex $1$ is one and leaving vertex $n+1$ is one, and there is conservation of fluid at the intermediate vertices. Formally, a \textbf{flow} $f$ of size one on $G$ is a function $f: E \rightarrow \R_{\geq 0}$ from the edge set $E$ of $G$ to the set of nonnegative real numbers such that 
  
  $$1=\sum_{e \in E, \ini(e)=1}f(e)= \sum_{e \in E,  \fin(e)=n+1}f(e),$$
  
  and for $2\leq i\leq n$
  
  $$\sum_{e \in E, \fin(e)=i}f(e)= \sum_{e \in E, \ini(e)=i}f(e).$$
  
  \medskip

  The \textbf{flow polytope} $\F_G$ associated to the graph $G$ is the set of all flows $f: E \rightarrow \R_{\geq 0}$ of size one. In their unpublished work \cite{p, S} Postnikov and Stanley  discovered a remarkable connection between the volume of the flow polytope and the Kostant partition function $K_G$. Namely, they proved that given a loopless  graph $G$ on the vertex set $[n+1]$, the normalized volume $\vol(\F_G)$ of the flow polytope associated to graph $G$ is 
\begin{equation} \label{eq:vol} \vol(\F_G)=K_G(0,  d_2, \ldots, d_{n}, -\sum_{i=2}^{n} d_i), \end{equation}
  where  $d_i=indeg_G(i)-1$ for $i \in \{2, \ldots, n\}$, and   $K_G(\v)$ denotes the \textbf{Kostant partition function}, which is the number of ways to write the vector $\v$ as a nonnegative linear combination of the positive type $A_n$ roots corresponding to the edges of $G$, without regard to order. To the edge  $(i, j)$, $i<j$, of $G$ corresponds  the positive type $A_n$ root $e_i-e_j$, where $e_i$ is the $i^{th}$ standard basis vector in $\mathbb{R}^{n+1}$.  

Along with Postnikov and Stanley \cite{p, S}, Baldoni and Vergne \cite{bv, BV} also studied type $A_n$  flow polytopes extensively with residue techniques. The author and Morales \cite{mm} worked on  flow polytopes of other types using combinatorial techniques.

\subsection{Triangulation of flow polytopes}
\label{sec:tri}

We show a systematic way of triangulating flow polytopes $\F_G$ associated to graphs $G$ with special source and sink vertices.  This process was used by \cite{p, S} in their unpublished work and appears in its most general form in \cite{mm}. A related process has been studied in detail in the context of root polytopes by the author \cite{m1, m2}.

 Let $G=([n+1], E)$  and let ${\tilde{G}}=([n+1]\cup \{s, t\}, \tilde{E})$, where $s$ is the smallest, $t$ is the biggest vertex of  $[n+1]\cup \{s, t\}$, and $\tilde{E}=E\cup \{(s, i), (i, t) | i \in [n+1]\}$.

   {\bf Algorithmic step:}    Given   a graph $G_0$ on the vertex set $[n+1]$ and   $(i, j), (j, k) \in E(G_0)$ for some $i<j<k$, let   $G_1$ and $G_2$ be graphs on the vertex set $[n+1]$ with edge sets
  \begin{eqnarray} \label{graphs1}
E(G_1)&=&E(G_0)\backslash \{(j, k)\} \cup \{(i, k)\}, \nonumber \\
E(G_2)&=&E(G_0)\backslash \{(i, j)\} \cup \{(i, k)\}. 
\end{eqnarray}

\begin{figure}
\begin{center}
\includegraphics[scale=.7]{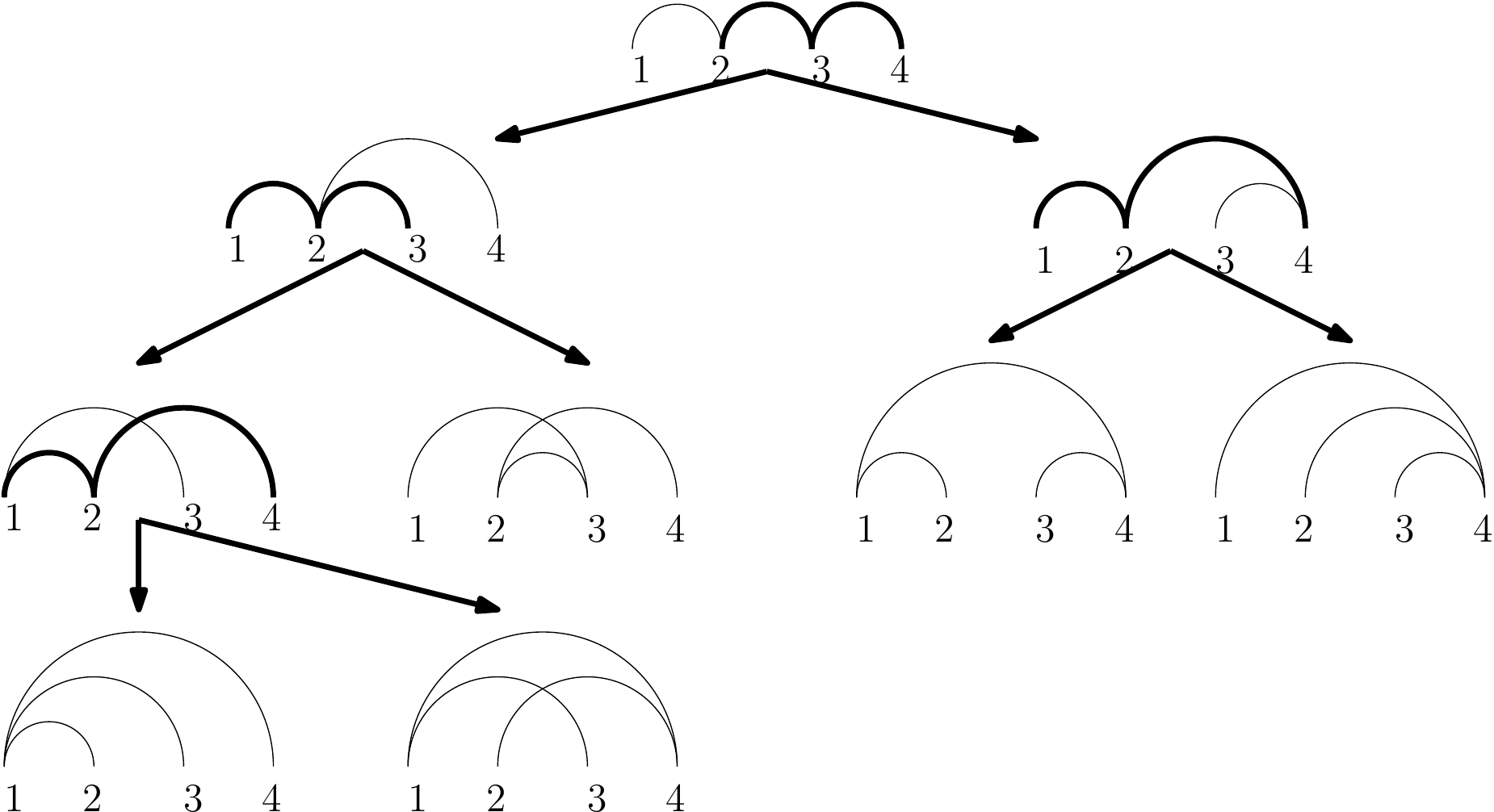}
 \caption{A reduction tree of $G=([4], \{(1,2),(2,3), (3,4)\})$. The edges on which the reductions are performed are in bold. }
 \label{fig}
 \end{center}
\end{figure}

    We say that $G_0$ \textbf{reduces} to $G_1$ and $G_2$ under the reduction rules defined by equations (\ref{graphs1}). A \textbf{reduction tree $\T(G)$} of a graph $G$ is defined as follows.  The root of $\T(G)$ is labeled by $G$. Each node $G_0$ in $\T(G)$  has two children, which depend on the choice of the edges  of  $G_0$ on which we perform the reduction. Namely, if the reduction  is performed on edges $(i, j), (j, k) \in E(G_0)$,  $i<j<k$, then the two children of  the node $G_0$   are labeled by the graphs   $G_1$ and  $G_2$   as described by equation (\ref{graphs1}). For an example of a reduction tree see Figure ~\ref{fig}. Note that the reduction tree of $G$ is not unique: it depends on the particular reductions we choose at each graph. However, the number of leaves of all reduction trees of $G$ is the same, and is the value of the normalized volume of $\F_{{\tilde{G}}}$ as Corollary \ref{norvol} states below.
   
 \begin{proposition} \cite{p, S} \label{red} Given a graph $G_0$ on the vertex set $[n+1]$ and   $(i, j), (j, k) \in E(G_0)$,   for some $i<j<k$, let $G_1$ and $G_2$ be as in equations (\ref{graphs1}). Then 
 $$\F_{{\tilde{G}_0}}=\F_{{\tilde{G}_1}} \bigcup \F_{{\tilde{G}_2}} \text{ and } \F_{{\tilde{G}_1}}^\circ \bigcap \F_{{\tilde{G}_2}}^\circ=\emptyset,$$

 \noindent where $\F_{{\tilde{G}_0}}$, $\F_{{\tilde{G}_1}}$, $\F_{{\tilde{G}_2}}$ are of the same dimension and  $\mathcal{P}^\circ$ denotes the interior of $\mathcal{P}$.
\end{proposition}

  For a  proof of Proposition \ref{red} see \cite{mm}.
 
 \begin{corollary} \label{norvol} The normalized volume of $\F_{{\tilde{G}}}$ is equal to the number of leaves in a reduction tree ~$\T(G)$.
 \end{corollary}

 \section{Volumes of flow polytopes and triangular arrays}
\label{sec:arr}
 
In this section we encode the leaves of a reduction tree $\T(G)$ in terms of triangular arrays. In light of Corollary \ref{norvol}, the calculation of the volume of  $\F_{{\tilde{G}}}$  is then a matter of enumerating these triangular arrays. We also show that the volume of $\F_{{\tilde{G}}}$ is equal to the constant term of  a formal Laurent series.

\begin{figure}[htbp] 
\begin{center} 
\includegraphics[scale=.5]{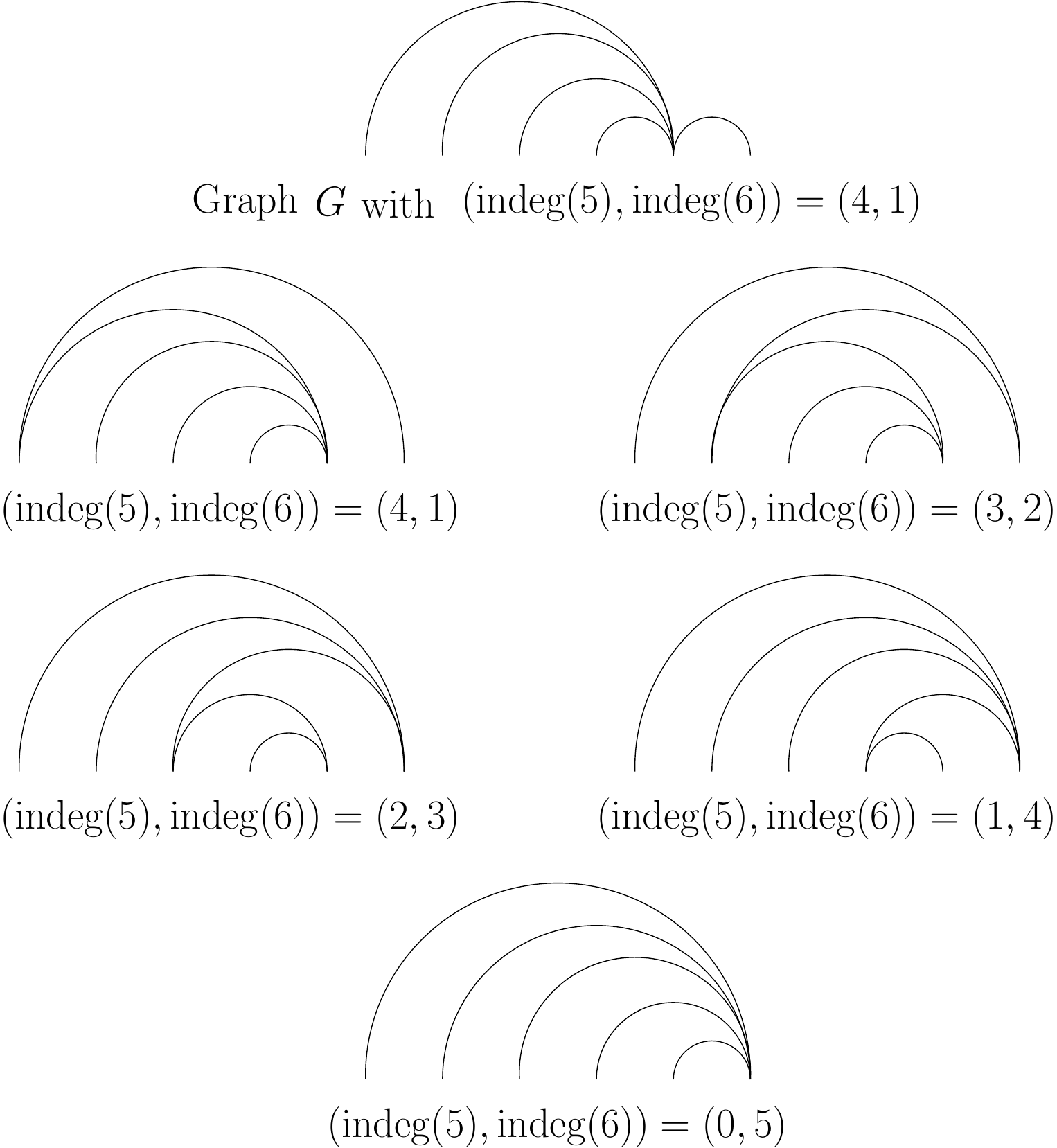} 
\caption{Graph $G$ and the graphs labeling the leaves of its reduction tree.} 
\label{indeg} 
\end{center} 
\end{figure}

 The basic idea is to  encode the number of incoming edges of the vertices. We illustrate this  idea on a graph $G$ on the vertex set $[a+2]$ with edges $(i, a+1)$, $i \in [a]$, and $(a+1, a+2)$. Figure \ref{indeg} shows the graphs labeling the leaves of the reduction tree of $G$ for $a=4$.   Note that the pair $(\i(a+1), \i(a+2))$ uniquely determines any leaf of the reduction tree and it takes each value $(a+1-i, i),$ $i \in [a+1]$, exactly once.

 The above simple example is the cornerstone of encoding the leaves of a reduction tree of any graph where the only possible multiple edges are of the form $(1, l)$. The case where multiple edges are allowed in general, can be dealt with similarly, and Lemma \ref{bas} is the cornerstone of that case. For simplicity, in this paper we mostly consider the case where possible multiple edges are of the form $(1, l)$. For the latter, the next lemma is key.

 \begin{lemma} \label{key}
 Given a graph $G$ and a distinguished vertex $v$ in it with $a$ incoming edges and one outgoing edge $(v, u)$, there is a way to perform all reductions possible which involve only edges incident to $v$, so that at the end of the process we obtain  graphs $G_i$, $i \in [a+1]$,  such that $(\i_{G_i}(v), \i_{G_i}(u))=(a+1-i, i).$ 
 \end{lemma}
 
  Given a graph $G$ on the vertex set $[n+1]$, associate to it its \textbf{indegree sequence}  $$\i(G)=(\i(2), \ldots, \i(n+1)).$$ For notational convenience we also define the  \textbf{$m$-indegree sequence} of $G$, $m \in \mathbb{Z}_{\geq 2}$,  to be $$\i^m(G)=(\i(2), \ldots, \i(m)),$$ where $$\i(k)=0 \text{ for } k>n+1.$$

 \medskip
 
We  recursively characterize the indegree sequences of the leaves of a particular reduction tree of $G$. We now describe our procedure for doing the reductions. 
 
Denote by  $I_i$ the set of incoming edges into vertex $i$ in the graph $G$. Denote by $V_i$ the vertices $k<i$ such that $(k, i)$ is an edge of  $G$. Denote by $G[i]$ the restriction of $G$ to the vertex set $[i]$. Denote by ${\rm InSeq}(\T(G))$ the multiset of indegree sequences of the graphs labeling the leaves of the reduction tree $\T(G)$, and by ${\rm InSeq}^m(\T(G))$ the multiset of $m$-indegree sequences of the graphs labeling the leaves of the reduction tree $\T(G)$. 

\medskip

\noindent \textbf{ALGORITHM 1.} Recursive description of ${\rm InSeq}(\T(G))$: 
\begin{itemize}

\item Construct the reduction tree $\T(G[2])$. Its only leaf is $G[2]$;  thus,  $${\rm InSeq}(\T(G[2]))= \{(\i(G[2]))\}=\{(|E(G[2])|)\}.$$

 \item Having constructed $\T(G[i])$, construct the  reduction tree $\T(G[i+1])$ from $\T(G[i])$ by appending the vertex $i+1$  and the edges $I_{i+1}$ to all graphs  in $\T(G[i])$, and then performing reductions at each vertex in $V_{i+1}$ on the graphs corresponding to the leaves of $\T(G[i])$, as described below. 
 \medskip
 
\item  Let  $V_{i+1}=\{i_1<i_2<\cdots<i_k\}$ and let $(s_2, \ldots, s_{n+1})$ be one of the sequences in ${\rm InSeq}^{n+1}(\T(G[i]))$. Applying Lemma \ref{key} a total of  $k$ times to the vertices  $i_1, \ldots, i_k$, in this order, we see that the leaves of  $\T(G[i+1])$ which are descendents of the graph which has $(n+1)$-indegree sequence  $(s_2, \ldots, s_{n+1})$ have $(n+1)$-indegree sequences

$$(s_2, \ldots, s_{n+1})+v^{i+1}[i_1]+\cdots+v^{i+1}[i_k],$$

where

\begin{align*} v^{i+1}[i_l] \in \{&(c_2, \ldots, c_{n+1})| c_i=0 \text{ for } i \in [n+1]\backslash \{i_l, i+1\}, c_{i_l}=s_{i_l}+1-s, \\ & c_{i+1}=s, \text{ for } s \in [s_{i_l}+1]\}, \text{ for } l \in [k].\end{align*}

 \end{itemize}

 The above procedure yields a complete description of the set ${\rm InSeq}(\T(G))$  using the sets $V_{i}$, $i \in [n+1]$. In order to calculate the volume of $\F_{{\tilde{G}}}$ we need to evaluate $|{\rm InSeq}(\T(G))|$, and  ALGORITHM 1 imples the following theorem.

 \begin{theorem}\label{thm:a}
 $$\vol(\F_{{\tilde{G}}})=|{\rm InSeq}(\T(G))|.$$
 Furthermore, $|{\rm InSeq}(\T(G))|$ can be characterized as follows.  Let $(a_{i,1}, \ldots, a_{i,i})$ be the indegree sequences of $\T(G[i+1])$, $ i \in [n]$, as defined above. Then, 
 \begin{itemize}
 
\item $a_{11}=|E(G[2])|$

\item $a_{i,j}\in \{0, \ldots, a_{i-1, j}\},$ if $(j+1, i+1)\in G$

\item $a_{i,j}= a_{i-1, j},$ if $(j+1, i+1) \not \in G$

\item $a_{i,i}=|E(G[i+1])|-\sum_{k=1}^{i-1}a_{i,k}$  
 \end{itemize}
 
 \end{theorem}

In the above, and below,  $a_{ij}$ denotes the variable for the indegree of vertex $j+1$ in the graph on $i+1$ vertices. 
  The variables $(a_{ij})_{i> j}$ can be arranged in a triangular array. For $$G=([5], \{(1,2), (1,3), (1,5), (2,4), (3,5)\})$$ the triangular array is

\hspace{2in} $ a_{4,1}\text{ } a_{31}\text{ }  a_{21}$

\hspace{2in} $ a_{4,2}\text{ } a_{32}$

\hspace{2in} $  a_{4,3}$

subject to the constraints

\hspace{2in} $ 0\leq a_{4,1}=a_{31}\leq a_{21}=1(=a_{11})$

\hspace{2in} $ 0\leq a_{42}\leq a_{32}=2-a_{11}(=a_{22})$

\hspace{2in} $ 0\leq a_{43}=3-a_{32}-a_{31}(=a_{33})$

\bigskip

 Theorem \ref{thm:a} states that $\vol(\F_{{\tilde{G}}})$ is equal to the number of triangular arrays with the constraints given above.
\bigskip

 \begin{theorem}\label{thm:b} The normalized volume 
 $\vol(\F_{{\tilde{G}}})$, where $G$ is a graph on the vertex set $[n+1]$,  is equal to the number of triangular arrays $(b_{i,j})_{i>j}$, $ j \in [n-1]$, $i \in \{j+1, \ldots, n\}$,   with the constraints

 $$b_{n,i} + b_{n-1, i} + \cdots +b_{i+1, i}  \leq  |E(G[i+1])|-|E(G[i])| +\sum_{k=1}^{i-1}b_{ik}, \text{ for all } i \in [n-1]$$
 
 and constraints  $b_{ji}=0$ if $(i+1, j+1)\not \in E(G)$.
 
 \end{theorem}

 \proof Let $$b_{ii}= |E(G[i+1])|-|E(G[i])| +\sum_{k=1}^{i-1}b_{ik}.$$

 Then the following is a bijection between the arrays  $(b_{i,j})_{i\geq j}$,  and $(a_{i,j})_{i\geq j}$, $ j \in [n-1]$, $i \in \{j, \ldots, n\}$, 
 
  $b_{ii}=a_{ii}=|E(G[i+1])|-\sum_{k=1}^{i-1}a_{ik}$
  
   $b_{ji}=a_{j-1, i}-a_{j,i}$ $i+1\leq j\leq n$
 
 Note that  $b_{ji}=0$ if $a_{j-1, i}=a_{j,i}$, that is if $(i+1, j+1)\not \in E(G)$.
  \qed
 
 \begin{theorem}  \label{b} The normalized volume 
 $\vol(\F_{{\tilde{G}}})$, where $G$ is a graph on the vertex set $[n+1]$, is equal to 
 
 $$CT_{x_n} \cdots CT_{x_1} \prod_{i=1}^n (1-x_i)^{-1}  \prod_{(i, n+1)\in E(G)} (1-x_i)^{-1}  \prod_{i=2}^n x_i^{-c_{i-1}} \prod_{(i, j)\in E(G)\mbox{ }: \mbox{ }j\leq n} (1-\frac{x_i}{x_j})^{-1}, $$
where $c_i= |E(G[i+1])|-|E(G[i])|$, and $CT_{x_i}$ mean the constant term in the expansion of the variable $x_i$.

\end{theorem}

\proof It suffices to show that the  number of triangular arrays  $(b_{ij})_{i> j}$ as described in  Theorem \ref{thm:b}, corresponding to graph $G$ is equal to  the constant term described above.

Given a  triangular array  $(b_{ij})_{i\geq j}$ corresponding to graph $G$ on the vertex set $[n+1]$, define $s_i$ by the equation

$$s_i+ \sum_{m=i+1}^n b_{m,i}   = |E(G[i+1])|-|E(G[i])| +\sum_{k=1}^{i-1}b_{ik}, \text{ for all } i \in [n-1].$$ Also, let $b_{i,0}=0$ for $i \in [n]$.

Then, the following term  of the above product is constant:

$$\prod_{i=1}^n x_i^{s_{i-1}} \prod_{(i, n+1)\in E(G)} x_i^{b_{n, i-1}}  \prod_{i=2}^n x_i^{-c_{i-1}} \prod_{(i, j)\in E(G)\mbox{ }: \mbox{ }j\leq n}  (\frac{x_i}{x_j})^{b_{j-1, i-1}},$$
since the power of  $x_{i+1}$, $i \geq 1$,  in it is $s_i+\sum_{m=i+1}^n b_{m,i}   -|E(G[i+1])|+|E(G[i])| -\sum_{k=1}^{i-1}b_{ik}=0,$ $i \in [n-1]$ and the power of $x_1$ is similarly $0$. Conversely, observe that any constant term is of the previous form, proving Theorem \ref{b}.
\qed
 
In the next sections we utilize the above theorems for calculating  volumes of flow polytopes.
 
 \section{Polytopes with volume  $\prod_{i=m+1}^{m+n-2}\frac{1}{2i+1}{{m+n+i} \choose {2i}}$}
 \label{sec:m}
 
In this section we define a family of polytopes $\mathcal{P}_{m,n}$, $m \in \mathbb{Z}_{\geq 0}$, $n \in \mathbb{Z}_{\geq 2}$, 
of which the Chan-Robbins-Yuen polytope is a special case, and prove that the normalized volume of $\mathcal{P}_{m,n}$ is $$\prod_{i=m+1}^{m+n-2}\frac{1}{2i+1}{{m+n+i} \choose {2i}}.$$ At the same time we also provide  a polytopal interpretation for  Conjecture 3 of Chan, Robbins and Yuen \cite[Conjecture 3]{cry} which was given only in terms of triangular arrays previously. Finally, we prove a  Kostant partition function identity as a corollary of our results. 
 
\begin{figure}[htbp] 
\begin{center} 
\includegraphics[scale=.65]{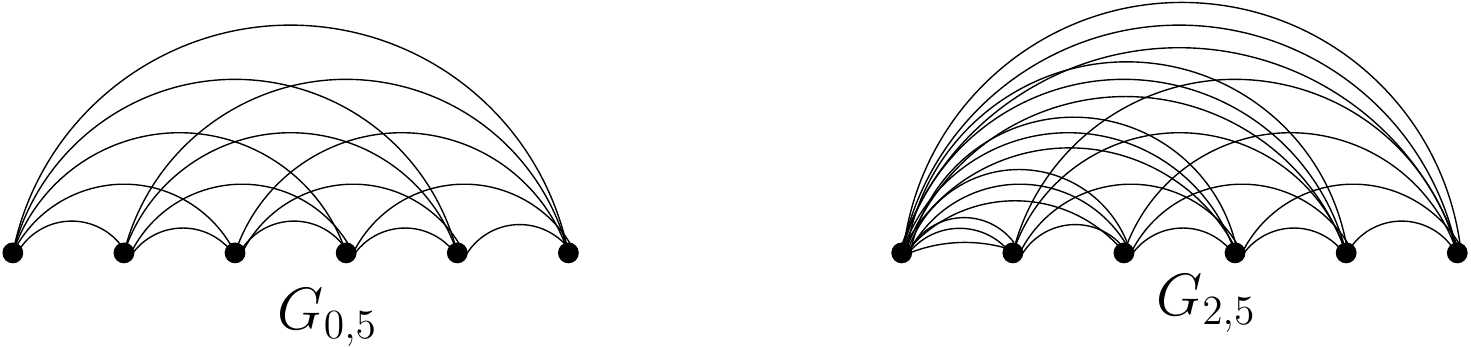} 
\caption{Graphs $G_{0,5}$  and $G_{2,5}$.} 
\label{G_0} 
\end{center} 
\end{figure}

  Consider the family of graphs $G_{m, n}$, $m \in \Z_{\geq 0}$, on the vertex set $[n+1]$ and with multiset of edges containing all edges of the complete graph, where the edges incident to $1$ have multiplicity $m+1$. See Figure \ref{G_0} for $G_{0,5}$ and $G_{2,5}$.

By Theorem \ref{thm:a},  $|{\rm InSeq}(\T(G_{m,n}))|$, $m \in \Z_{\geq 0}$, is enumerated  by the following triangular arrays.  Let $(a_{i1}, \ldots, a_{ii})$ be the indegree sequences of $\T(G_{m,n}[i+1])$. Then, 
 \begin{itemize}
 
\item $a_{11}=m+1$

\item $a_{ij}\in \{0, \ldots, a_{i-1, j}\},$ for $2\leq i$, $1\leq j<i$

\item $a_{ii}={{i+1} \choose 2}+im-\sum_{k=1}^{i-1}a_{ik}$  
 \end{itemize}
 
We can  write the conditions on the $a_{ij}$, $i, j \in [n]$, $j<i$,  in a triangular array. Namely, let the $i^{th}$ row of the triangular array be

 $0\leq a_{n,i} \leq a_{n-1, i} \leq \cdots \leq a_{i+1, i}  \leq  { i+1 \choose 2} +ib-\sum_{k=1}^{i-1}a_{ik}$
 \medskip

  The triangular array for $n=5$ looks like:
  
\hspace{2in} $0\leq a_{5,1} \leq a_{4, 1}\leq a_{3, 1}  \leq a_{21}\leq { 2 \choose 2} +m$

\hspace{2in} $0\leq a_{5,2} \leq a_{4, 2}\leq a_{3, 2} \leq { 3 \choose 2} +2m-a_{21}$

\hspace{2in} $0\leq a_{5,3} \leq a_{4, 3} \leq { 4 \choose 2} +3m-a_{31}-a_{32}$

\hspace{2in}  $0\leq a_{5,4}  \leq { 5 \choose 2} +4m-a_{41}-a_{42}-a_{43}$
 \bigskip
 
 The following is a special case of Theorem \ref{thm:b}.
 \begin{theorem} \label{thm:bb}
 There is a bijection between triangular arrays  whose $i^{th}$ row  $a_{n,i},  a_{n-1, i}, \cdots,  a_{i+1, i}$ with $a_{ij}\geq 0$,
 satisfies

 $$a_{n,i} \leq a_{n-1, i} \leq \cdots \leq a_{i+1, i}  \leq  { i+1 \choose 2} +im-\sum_{k=1}^{i-1}a_{ik}$$
 
\noindent and triangular arrays  whose $i^{th}$ row  $b_{n,i},  b_{n-1, i}, \cdots,  b_{i+1, i}$, with $b_{ij}\geq 0$,
 satisfies

$$b_{n,i} + b_{n-1, i} + \cdots +b_{i+1, i}  \leq  m+i +\sum_{k=1}^{i-1}b_{ik}.$$
 
 \end{theorem}
 
 Using Theorem \ref{thm:bb} we arrive at the main result of this section, which is a formula for the volume of $\F_{{\tilde{G}_{m,n}}}$ as well as the geometrical realization of \cite[Conjecture 3]{cry} as the proof  of Theorem \ref{prod} shows:
  
 \begin{theorem} \label{prod} The normalized  volume of the flow polytope $\F_{{\tilde{G}_{m,n}}}$ is
 $$\prod_{i=m+1}^{m+n-2}\frac{1}{2i+1}{{m+n+i} \choose {2i}}$$
 
 \end{theorem}
 
 \proof  By Theorems \ref{thm:a} and  \ref{thm:bb}, the normalized  volume of the flow polytope $\F_{{\tilde{G}_{m,n}}}$ is equal to the number of triangular array  $(b_{ij})_{j<i}$ as described in Theorem \ref{thm:bb}.
The triangular array  $(b_{ij})_{j<i}$ agrees with the set of triangular arrays $\mathcal{A}_{m+n+1, m+1}$ described in \cite[Conjecture 3]{cry} upon appending $m+1$ rows of $0'$s on top. \cite[Conjecture 3]{cry}  states that the number of such arrays is $\prod_{i=m+1}^{m+n-2}\frac{1}{2i+1}{{m+n+i} \choose {2i}}$. Zeilberger proved this in \cite{z}. \qed

 \medskip
 
 Theorem \ref{prod} implies the Kostant partition function identity stated in Corollary \ref{kir}. This identity has been previously observed by A.N. Kirillov in a somewhat different form and is stated in his paper \cite[p. 80 ]{kir} without proof.  
 
 \begin{corollary} \label{kir}
 \begin{equation} \label{eq:kost}K_{K_{n+1}}(m+1, m+2, \ldots, m+n, -nm-{n \choose 2})=\prod_{i=m+1}^{m+n-2}\frac{1}{2i+1}{{m+n+i} \choose {2i}}.\end{equation}
 \end{corollary}
 
 \proof By \eqref{eq:vol}, Theorem \ref{prod} implies
 
$$K_{{\tilde{G}_{m,n}}}(0,  0, m+1, m+2, \ldots, m+n, -nm-{n \choose 2})=\prod_{i=m+1}^{m+n-2}\frac{1}{2i+1}{{m+n+i} \choose {2i}},$$ which  can be easily seen to be equivalent to  \eqref{eq:kost}.
 \qed

 Note that the right hand side of \eqref{eq:kost} can also be written as $$\prod_{p=1}^{n-2} C_p \prod_{1\leq i<j\leq n-1} \frac{2(m+1)+i+j-1}{i+j-1}$$ which is the way it appears in \cite[p. 80]{kir}.

   \section{A collection of  polytopes whose volumes sum to $N(n,k)\times \prod_{i=1}^{n-1}C_i$}
\label{sec:conj2}

In this section we construct for each $n, k \in \Z$, $k \leq n$,  a collection of polytopes $\{P^{n,k}_i\}_{i \in C_{n,k}}$ whose volumes sum to $$N(n,k)\times \prod_{i=1}^{n-1}C_i,$$ where $$N(n, k)=\frac{1}{n}{n \choose k}{n \choose {k-1}}$$ is the Narayana number and $C_i$ denotes the $i$th Catalan number.  While curious on its own, we do the above in order to  provide a polytopal interpretation of \cite[Conjecture 2]{cry}.

The collections $\{P^{n,k}_i\}_{i \in C_{n,k}}$ are such that for fixed $n$, the polytopes $P^{n,k}_i$,  $k=1,2,\ldots, n$, $i \in C_{n,k}$, are interior disjoint and their union is $CRY_{n+2}$.

 It follows form the definition of  $CRY_{n+2}$ that it can be thought of as a flow polytope  of the graph ${\tilde{K}_{n+1}}$. In particular  $$\vol(\F_{{\tilde{K}_{n+1}}})=\prod_{i=1}^{n}C_i=\sum_{k=1}^n  N(n,k)\times \prod_{i=1}^{n-1}C_i.$$
 
 Let the collection $\{G^{n,k}_i\}_{i \in C_{n,k}}$ of graphs be all graphs on the vertex set $[n+1]$ consisting of ${{n+1} \choose 2}$ edges, such that
 
 \begin{itemize}
 \item  the edge $(1, n+1)$ has multiplicity $n$
 
 \item all edges other than $(1, n+1)$ have multiplicity zero or one
 
 \item  if edge $(1, l)$, $2\leq l\leq n$, is in the graph, then edge $(l, n+1)$ is not
 
 \item there are $k-1$ edges of the form $(1, l)$, $2\leq l\leq n$, in the graph
 
 \end{itemize}
 
\begin{figure}[htbp] 
\begin{center} 
\includegraphics[scale=.7]{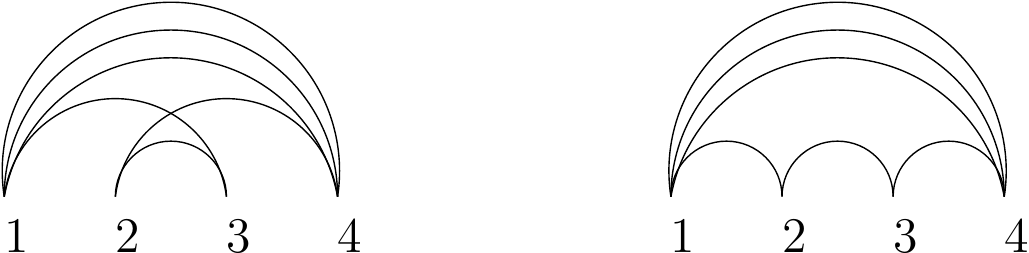} 
\caption{Graphs  $\{G^{3,2}_i\}_{i \in C_{3,2}}.$} 
\label{t3} 
\end{center} 
\end{figure}

       See Figure \ref{t3} for an illustration of the graphs  $\{G^{3,2}_i\}_{i \in C_{3,2}}.$

 Define $$P^{n,k}_i=\F_{{\tilde{G}^{n,k}_i}}, \text{ for } {k\leq n, i \in C_{n,k}}.$$ We set up the definition of $\{P^{n,k}_i\}_{i \in C_{n,k}}$ exactly so that for fixed $n$, the polytopes $P^{n,k}_i$,  $k=1,2,\ldots, n$, $i \in C_{n,k}$, are interior disjoint and their union is $\F_{{\tilde{K}_{n+1}}}$. 
 
 \begin{proposition} \label{sub}
 For fixed $n$, the polytopes $P^{n,k}_i$,  $k=1,2,\ldots, n$, $i \in C_{n,k}$, are interior disjoint and their union is $\F_{{\tilde{K}_{n+1}}}$. 
 \end{proposition}
 
 \proof
 Note that the graphs  $G^{n,k}_i$,  $k=1,2,\ldots, n$, $i \in C_{n,k}$, can be obtained from $K_{n+1}$ by applying the reduction rules as specified in \eqref{graphs1} on the edge pairs $(1,l)$, $(l, n+1)$, for each $l \in \{2,3,\ldots, n\}$. Thus, Proposition \ref{red} implies Proposition \ref{sub}.
 \qed
 
 \medskip
 
Next, we apply Theorem \ref{thm:b} to encode the normalized volumes of the polytopes $P^{n,k}_i=\F_{{\tilde{G}^{n,k}_i}},$  for  ${k\leq n, i \in C_{n,k}},$ by triangular arrays. 

\begin{theorem} \label{Pnk}
Fix $n, k \in \Z$ such that $1\leq k \leq n$. Then, the sum of the normalized volumes of the polytopes $P^{n,k}_i, i \in C_{n,k}$, is equal to the number of triangular arrays $(b_{i,j})_{i>j}$, $ j \in [n-1]$, $i \in \{j+1, \ldots, n\}$,   with the constraints

 $$b_{n,i} + b_{n-1, i} + \cdots +b_{i+1, i}  \leq  n_i +\sum_{k=1}^{i-1}b_{ik}, \text{ for all } i \in [n-1],$$
 such that  $k-1$ of the variables $b_{n,j}$, $1\leq j\leq n-1$  are required to be $0$ a priori,  and $n_i=i$ if $b_{n,i}$ is among these $k$ variables, and $n_i=i-1$ otherwise.

\end{theorem}

\proof The proof of Theorem \ref{Pnk} is a straightforward application of Theorem \ref{thm:b} for each  polytope $P^{n,k}_i=\F_{{\tilde{G}^{n,k}_i}},  i \in C_{n,k},$ 
 involved. 
\qed

\begin{lemma} \label{lem:tild} The
triangular arrays $(b_{i,j})_{i>j}$, $ j \in [n-1]$, $i \in \{j+1, \ldots, n\}$,   with the constraints 
 
 $$b_{n,i} + b_{n-1, i} + \cdots +b_{i+1, i}  \leq  n_i +\sum_{k=1}^{i-1}b_{ik}, \text{ for all } i \in [n-1]$$
 such that  $k-1$ of the variables $b_{n,j}$, $1\leq j\leq n-1$  are required to be $0$ a priori,  and $n_i=i$ if $b_{n,i}$ is among these $k$ variables, and $n_i=i-1$ otherwise,  are in bijection with triangular arrays  $(\tilde{b}_{i,j})_{i>j}$, $ j \in [n-1]$, $i \in \{j+1, \ldots, n\}$,   with the constraints

 $$\tilde{b}_{n,i} + \tilde{b}_{n-1, i} + \cdots +\tilde{b}_{i+1, i}  \leq  i +\sum_{k=1}^{i-1}\tilde{b}_{ik}, \text{ for all } i \in [n-1]$$
 such that  exactly $k-1$ of the inequalities hold at equality.
\end{lemma}

\proof Given a triangular arrays $(b_{i,j})_{i>j}$, $ j \in [n-1]$, $i \in \{j+1, \ldots, n\}$, of the first kind with $J \subset [n-1]$, $|J|=k-1$, such that the variables $b_{n,j}$, $j \in J$, are required to be  $0$, let $(\tilde{b}_{i,j})_{i>j}$, $ j \in [n-1]$, $i \in \{j+1, \ldots, n\}$, be the array obtained by setting $\tilde{b}_{i,j}=b_{i,j}$ for $(i, j)\neq (n,l)$ for $l \in J$. Finally, let $\tilde{b}_{n,j}=j-\sum_{j+1}^{n-1}\tilde{b}_{m,j}$, for $j \in J$. Note that this way we obtained a triangular array $(\tilde{b}_{i,j})_{i>j}$ for which exactly $k-1$ inequalities hold at equality, namely, the rows indexed by $j \in J$. 
\qed

\begin{theorem} \label{thm:nar}
Fix $n, k \in \Z$ such that $1\leq k \leq n$. Then, the sum of the normalized volumes of the polytopes $P^{n,k}_i, i \in C_{n,k}$, is equal to $$N(n,k)\times \prod_{i=1}^{n-1}C_i.$$
\end{theorem}

\proof Conjecture 2 in the paper of Chan-Robbins-Yuen \cite{cry}, proved by Zeilberger \cite{z} enumerates the arrays $(\tilde{b}_{i,j})_{i>j}$, $ j \in [n-1]$, $i \in \{j+1, \ldots, n\}$ as defined in Lemma \ref{lem:tild}. Thus, Theorem \ref{thm:nar} follows from Theorem \ref{Pnk}, Lemma \ref{lem:tild} together with \cite[Conjecture 2]{cry} and its proof in \cite{z}.
\qed

\medskip

Theorem \ref{thm:nar} thus establishes a polytopal realization of \cite[Conjecture 2]{cry}. 

   \section{Constructing polytopes with volume  $ \frac{1}{(r-1)n+1} {{rn} \choose n}$}
\label{sec:constr}

In this section we give an example of how to use the triangular arrays to construct polytopes with combinatorial volumes. There is no prescribed formula for this, rather we use the triangular array encoding as stated in Theorem \ref{thm:b}  as  a  heuristic.

Say we would like to construct a family of polytopes with volumes equal the number of $(r+2)$-ary trees with $n+1$ internal nodes. If we are to use Theorem \ref{thm:b}, this amounts to finding a triangular array $(b_{i,j})_{i>j}$ with constraints as described in Theorem \ref{thm:b} and which bijects with $(r+2)$-ary trees with $n+1$ internal nodes. The following proposition takes us close to the solution.

\begin{proposition} \label{r-ary}
There is a bijection between  nonnegative integer sequences $(b_{i+1,i})_{i=1}^n$ with constraints $b_{2,1} \leq r+1$ and $b_{i+1,i}\leq r+1+b_{i, i-1}$ for $2\leq i\leq n$ and $(r+2)$-ary trees with $n+1$ internal nodes.
\end{proposition}

\begin{figure}[htbp] 
\begin{center} 
\includegraphics[scale=.9]{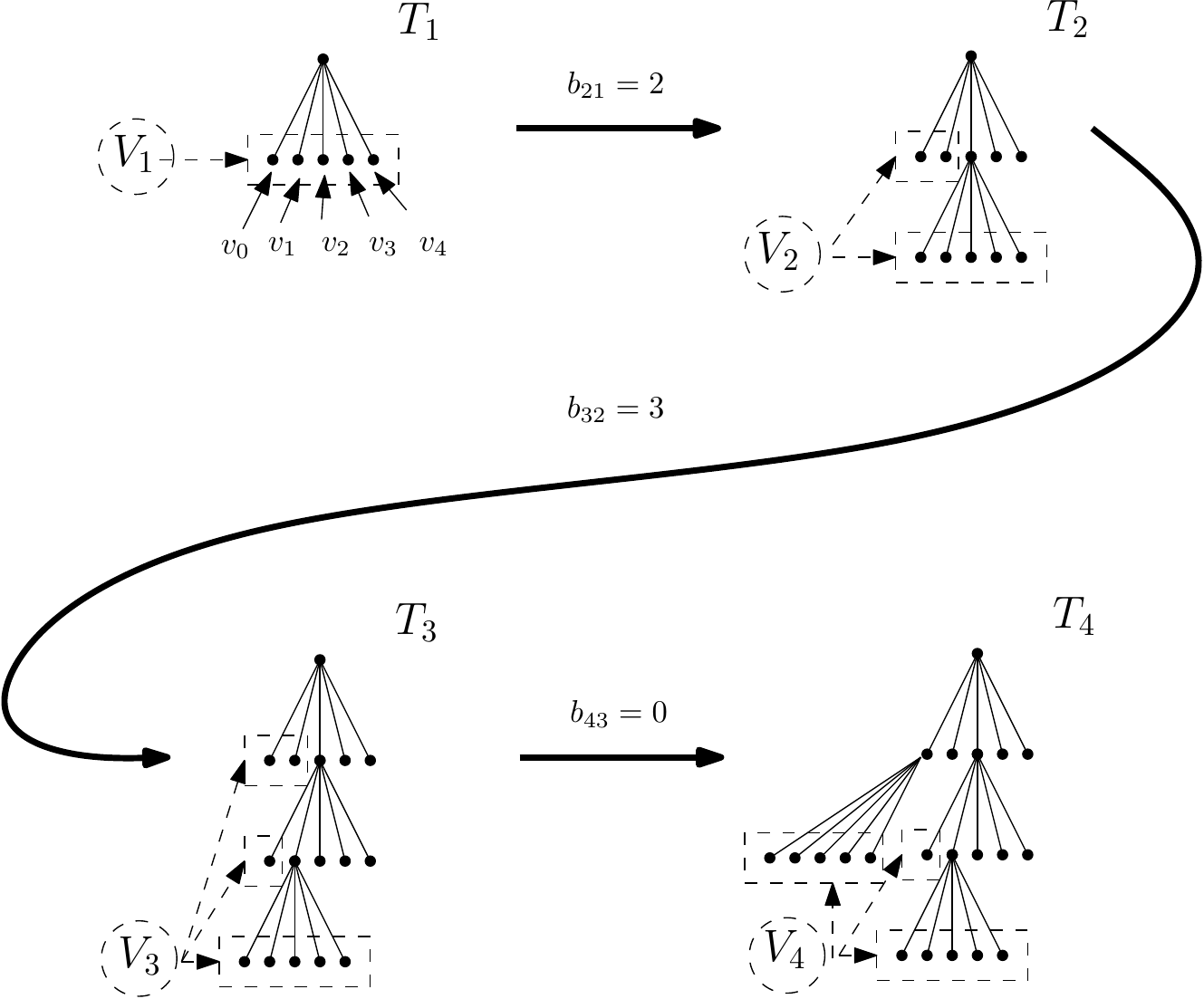} 
\caption{The definitions from the proof of Proposition \ref{r-ary}, as well as the construction for $n=3$, $r=3$, $(b_{21}, b_{32}, b_{43})=(2,3,0)$ are illustrated.} 
\label{t1} 
\end{center} 
\end{figure}
\proof Given a nonnegative integer sequences $(b_{i+1,i})_{i=1}^n$ with constraints $b_{2,1} \leq r+1$ and $b_{i+1,i}\leq r+1+b_{i, i-1}$ for $2\leq i\leq n$ construct a sequence of rooted ordered $(r+2)$-ary trees $T_1, \ldots, T_{n+1}$ as follows.  Having constructed $T_i$, $i \in [n]$, we will let $T_{i+1}$ be the tree obtained from  $T_i$  by marking  a leaf of it in a manner described later and adding $r+2$ children to the marked vertex.  Let $T_1$ be the  $(r+2)$-ary tree with $r+3$ vertices consisting of the  root $v$ and its children $v_0, v_1, \ldots, v_{r+1}$  ordered from left to right. This ordered set of vertices of cardinality $r+2$ is denoted by $V_1$.  Let $V_i$ be the set of vertices in $T_i$ which are leaves, and such that all vertices which are children of their parents and are to the left of them are also leaves. See Figure \ref{t1} for an example. Once the vertices of $V_i$ are ordered, and one of its vertices  $w$ is marked, we obtain $V_{i+1}$ by deleting the vertices in $V_i$ which are children of the parent of $w$ and are equal to $w$ or are to the right of $w$ and adding the $r+2$ children of $w$. The ordering on the vertices which were also in $V_i$  is inherited from $V_i$, and the $r+2$ new vertices are the last $r+2$ ones and are ordered according to their order from left to right. The marking of a vertex in $V_i$, $i \in [n]$ is done as follows.  Mark the $(b_{i+1, i}+1)st$ vertex of $V_i$.
By construction $|V_1|=r+2$ and $|V_i|=r+2+b_{i, i-1}$ for $2\leq i\leq n$. Thus, the constraints $b_{2,1} \leq r+1$ and $b_{i+1,i}\leq r+1+b_{i, i-1}$ for $2\leq i\leq n$ can also be written as $b_{2,1}+1 \leq |V_1|$ and $b_{i+1,i}+1\leq |V_i|$ for $2\leq i\leq n$. Since the above procedure is clearly invertible, Proposition \ref{r-ary} follows. 
\qed

\begin{proposition} \label{ary}
There is a bijection between  nonnegative integer sequences $(b_{i+1,i})_{i=1}^n$ with constraints $b_{2,1} \leq r+1$ and $b_{i+1,i}\leq r+1+b_{i, i-1}$ for $2\leq i\leq n$ and the triangular array $(b_{i,j})_{i>j}$ with constraints as described in Theorem \ref{thm:b} arising from the graph $G=([n+2], E),$ where the edge set consists of $r+1$ edges $(1,2)$, $r$ edges $(1, i)$, for $i \in \{3, 4, \ldots, n+2\}$, and one edge $(i-1, i)$, for $i \in \{3, 4, \ldots, n+2\}$. 
\end{proposition}
 
 \proof The bijection is  trivial since for the graph $G$ described above the triangular array $(b_{i,j})_{i>j}$ is such that  $b_{i,j}=0$ unless $(i, j)=(j+1,j)$, $j \in [n]$,  and we have the  constraints $b_{2,1} \leq r+1$ and $b_{i+1,i}\leq r+1+b_{i, i-1}$ 
for $2\leq i\leq n$. \qed

 \begin{corollary}
\begin{equation} \vol \F_{{\tilde{G}}}= \frac{1}{(r+1)(n+1)+1} {{(r+2)(n+1)} \choose {n+1}},
\end{equation} where $G$ is as in Proposition \ref{ary} above. \end{corollary} 
 
 \proof Immediate since the number of $(r+2)$-ary trees on $n+1$ internal vertices is given by $\frac{1}{(r+1)(n+1)n+1} {{(r+2)(n+1)} \choose {n+1}}$.\qed
 
 \begin{corollary}
\begin{equation} K_{P_{n+2}}(r+1, r+1, \ldots, r+1, -(n+1)(r+1))= \frac{1}{(r+1)(n+1)+1} {{(r+2)(n+1)} \choose {n+1}},
\end{equation} where $P_{n+2}=([n+2], \{(i, i+1) \mid i \in [n+1]\})$. 
 \end{corollary}

\proof Use equation \eqref{eq:vol}. \qed
 
 \section{The generalized Catalanotope}
\label{sec:cat}

The purpose of this section is twofold. On one hand, we study the generalization of the Catalonotope, introduced by  Postnikov and Stanley \cite{p, S} in unpublished work. On the other hand we present Lemma \ref{bas} which can be used to encode the triangulations of flow polytopes corresponding to graphs with multiple edges similarly to how we did it above in the case where multiple edges were of the form $(1, l)$.

Postnikov and Stanley \cite{p, S} studied a root polytope, which they called the Catalonotope in their unpublished work. We  consider the flow polytope $\F_{\tilde{P}_{n+1}}$ corresponding to the path graph $P_{n+1}=([n+1], \{(i, i+1)\mid i \in [n]\})$ a \textbf{Catalonotope}, since the triangulations of this polytope and the root polytope Catalonotope are in bijection. 
A natural generalization of the Catalonotope is then the flow polytope arising a graph  whose set of edges is  the same as that of  $P_{n+1}$, but the multiplicities are not necessarily all $1$.

Denote by $f(c_1, \ldots, c_n)$ the volume of $\F_{{\tilde{G}}_{c_1,\ldots, c_n}}$, where $G_{c_1,\ldots, c_n}$ is a graph on vertex set $[n+1]$ and contains $c_i$ edges $(i, i+1)$, for $i \in [n]$. To calculate the volume  $f(c_1, \ldots, c_n)$, as well as to be able to encode the triangulations of a flow polytope corresponding to a graph with multiple edges, we need the following lemma.

\begin{lemma} \label{bas}
Let $G$ be  the graph on the vertex set $[3]$ with $c_1$ edges  $(1, 2)$ and $c_2$ edges  $(2,3)$. Then the multiset $\is(G)$ consists of ${{c_1+c_2-1-i} \choose {c_2-1}}$ copies of $( i, c_1+c_2-i )$, for $i\in \{0, 1, 2, \ldots,c_1\}$.  
\end{lemma}

The proof of Lemma \ref{bas} is routine and is left to the reader.

Then Lemma \ref{bas} implies that 

$$f(c_1, c_2, \ldots, c_n)=\sum_{i=0}^{c_1}{{c_1+c_2-1-i} \choose {c_2-1}}f(c_1+c_2-i, c_3, c_4, \ldots, c_n)$$

Continually applying Lemma \ref{bas} we arrive to the formula

\eq \label{vol} f(c_1, c_2, \ldots, c_n)=\sum_{(k_2, \ldots, k_n)\in K({c_1,\ldots, c_n})} (\prod_{i=2}^{n-1}  {{c_i-1+k_i} \choose {c_i-1}}) {{c_n+k_{n}} \choose {c_n}},\en

where $$K({c_1,\ldots, c_n})=\{(k_2, \ldots, k_{n})| 0\leq k_2\leq c_1, 0\leq k_{i}\leq c_{i-1}+k_{i-1}, \text{ for }i \in \{3, \ldots, n-1\}, k_n=c_{n-1}+k_{n-1}\}.$$

We can rewrite equation (\ref{vol}) as stated in the theorem below.

\begin{theorem} 
\eq
\lb{vform1}
\vol \F_{{{\tilde{G}}_{c_1,\ldots, c_n}}} = \sum_{(k_2, \ldots, k_n)\in K({c_1,\ldots, c_n})} (\prod_{i=2}^{n-1} {(c_i) ^ {(k_i)} \over k_i!}) {(c_n+1) ^ {(k_n)} \over k_n!},
\en
where 

$$(a) ^ {(b)}:=\prod_{j=0}^{b-1} (a+j),$$ and $K({c_1,\ldots, c_n})$ is as defined above.

\end{theorem}

\begin{lemma} \label{cat} $$|K(\mathbf{1}_n)|=C_n,$$ where $\mathbf{1}_n=(1, \ldots, 1) \in \mathbb{Z}^n$, $C_n=\frac{1}{n+1}{ 2n \choose n}$ denotes the $n$th Catalan number.
\end{lemma}

\proof Let $f(n)=|K(\mathbf{1}_n)|$. 
A sequence $(k_2, \ldots, k_{n})\in K(\mathbf{1}_n)$ can be though of as specifying the steps of a walk on the nonnegative $x$-axis. Namely, start at the origin $(0,0)$ and on the $i$th step step to $(0, k_{i+1})$, $i \in [n-1]$.  Write $f(n)=f_1(n)+f_2(n)$, where $f_1(n)$ is the number of sequences $(k_2, \ldots, k_{n})\in K(\mathbf{1}_n)$ with $k_2=0$ and  $f_2(n)$ is the number of sequences $(k_2, \ldots, k_{n})\in K(\mathbf{1}_n)$ with $k_2=1$. Note that by definition  $f_1(n)=f(n-1)$. Consider a  sequence $s=(k_2, \ldots, k_{n})\in K(\mathbf{1})$ with $k_2=1$. Let the $(i+1)$st step be the first one where we are back to origin; that is $k_{i+2}=0$ and $k_j>0$, $i+2>j$, $i>0$. This means that the sequence $(k_3-1, k_4-1, \ldots k_{i+1}-1)$ is in  $K(\mathbf{1}_{i})$ and $(k_{i+3}, k_{i+2}, \ldots k_{n})$ is in  $K(\mathbf{1}_{n-i-1})$. Thus, 

$$f(n)=\sum_{i=0}^{n-1}f(i)f(n-i-1),$$ where we set $f(0)=1$. Since the Catalan numbers $C_n$ satisfy the same relation and the starting values are the same, it follows that $f(n)=C_n$.
\qed

\begin{corollary} \cite{p, S}

$$\vol(\F_{{\tilde{P}_{n+1}}})=C_n,$$ where $P_{n+1}=([n+1], \{(i, i+1) \mid i \in [n]\})$.
\end{corollary}

\proof
Follows immediately from the above, since $f(\mathbf{1}_n)=|K(\mathbf{1}_n)|=C_n$.
\qed

\medskip

Stanley \cite{sta} pointed out the vague resemblance of \eqref{vform1} to the volume of the Stanley-Pitman polytope \cite{sp}:

\eq
\lb{vform}
V(x_1, \ldots, x_n) = \sum_{\skk \in K_\dd} \,
\prod_{i=1}^\dd {x_i ^ {k_i} \over k_i!}
\en
where 
\eq
\lb{kdef}
K_\dd:= \{ \kk \in \nnints^\dd: \sum_{i=1}^j k_i \ge j \mbox{ \rm for all }
1 \le i \le \dd-1  \mbox{ \rm and }
\sum_{i=1}^\dd k_i = \dd \}
\nonumber \en
with $\nnints := \{0,1,2, \ldots \}$.

\bigskip

 Is there a connection between the two polytopes that makes  this resemblance  precise and endows it with  a meaning?

\section*{Acknowledgement} I would like to thank Richard Stanley, Alex Postnikov and Alejandro Morales for numerous conversation on the topic of flow polytopes and Kostant partition functions over the course of the past several years.

\end{document}